\DeclareMathAlphabet{\Bbl}{U}{msb}{m}{n}
\DeclareMathAlphabet{\goth}{U}{euf}{m}{n}

\documentclass{article}

 \newtheorem{thm}{Theorem}
 \newtheorem{cor}[thm]{Corollary}
 \newtheorem{lem}[thm]{Lemma}
 \newtheorem{prop}[thm]{Proposition}
 \newtheorem{defn}[thm]{Definition}
 \newtheorem{rem}[thm]{Remark}

\newcommand{\noi}{\noindent}
\newcommand{\vol}{\hbox{vol}}
\newcommand{\tr}{\hbox{tr}}
\newcommand{\ul}{\underline}
\newcommand{\ra}{\rightarrow}
\newcommand{\ol}{\overline}
\newcommand{\eps}{\epsilon}
\newcommand{\be}{\begin{enumerate}}
\newcommand{\ee}{\end{enumerate}}

\begin{document}

\title{Immobilization of convex bodies in $\Bbl{R}^n$}

\author{Anthony David Gilbert and Saul Hannington Nsubuga}



\date{October, 2018}
\maketitle

\begin{abstract}
\noi We extend to arbitrary finite $n$ the notion of immobilization of a convex body $O$ in $\Bbl{R}^n$ by a finite set of points $\mathcal{P}$ in the boundary of $O$. Because of its importance for this problem, necessary and sufficient conditions are found for the immobilization of an $n$-simplex. A fairly complete geometric description of these conditions is given: as $n$ increases from $n = 2$, some qualitative difference  in the nature of the sets $\mathcal{P}$ emerges.
\end{abstract}


\section{Introduction} \label{intro}
Immobilization problems were introduced by \cite{KUP} and were motivated by the need to understand the best position a machine (robot hand) can grasp an object, \cite{XM}.  There is now an extensive literature in Robotics journals on immobilization, for example \cite{CSU}, \cite{BMU}, \cite{mark2}, \cite{BFMM}, \cite{rimon1}, \cite{rimon2}, \cite{rimbud}. In \cite{CSU} it was proved that four points suffice to immoblize any plane body and 2$d$ points to immobilize a $d-$dimensional polytope. \cite{BMU} answered Kuperberg's conjecture in the affirmative proving that apart from a circular disk, any plane convex object with smooth boundary could be immobilized with three points. It was \cite{BFMM} who first brought out both the geometrical and algebraic aspects of  immobilization by treating the case of a tetrahedron. Their results have recently been used in \cite{BM} to prove a necessary condition on rotors in tetrahedra. An overview of the classical results on immobilization focusing on analysis, existence and synthesis is given in \cite{strappen}.

\medskip
\noi For obvious practical reasons, much of the literature focuses on the problem in $\Bbl{R}^3$, and for the actual grasping of objects, robots take full advantage of the effects of friction. The aim of this article is quite different. Following \cite{BFMM}, we study the purely geometric problem of grasping a smooth convex body $O$ at a finite set of points $\mathcal{P}$ in its boundary and seek conditions so that $O$ is completely immobilized. \cite{BFMM} examined this problem in $\Bbl{R}^3$ and found a set of necessary and sufficient conditions for immobilizing a $3$-simplex, i.e a tetrahedron; for reasons explained below, simplices play a particularly important role in the analysis. In this article, we reproduce their results, but generalise the problem to $\Bbl{R}^n$. We follow the approach of \cite{BFMM} by recasting the immobilization problem of the simplex $\Delta$,
 now in $\Bbl{R}^n,$ as an extremal problem, but thereafter our techniques are very different. Solution of the extremal problem leads to necessary and sufficient conditions for immobilization of $\Delta$ in terms of a matrix $A$ which encodes the geometry of $\Delta$ and the points $\mathcal{P};$ these conditions, which are algebraic in character, agree with the results of \cite{CSU} and \cite{BFMM} for $n = 2, 3$ respectively and reveal an interesting difference in behaviour between $n \leq 3$ and $n > 3.$ We conclude the article by interpreting the conditions on $A$ geometrically, and thus show that there is a particular set $\mathcal{P}$ which is optimal for immobilization in a sense described below. Some of the detailed proofs in the article have been relegated to the appendices.

\section{Immobilization of Convex Bodies and Simplices} \label{westart}
\noindent Let $O$ be a convex body in $\Bbl{R}^n$. We shall call a point $p$ in the boundary of $O$ a contact point, and a set $\mathcal{P}$ of points in the boundary of $O$ will be called a contact set. A contact set $\mathcal{P}$ is said to immobilize $O$ if, whenever $O$ is held fixed, any rigid motion of the points of $\mathcal{P}$ causes one or more of its points to penetrate into the interior of $O$; equivalently, penetration by at least one point $p \in \mathcal{P}$ occurs whenever $O$ is moved and the set $\mathcal{P}$ is held fixed; in that case, we will call $\mathcal{P}$ an immobilizing contact set.

\noindent For a smooth convex body in $\Bbl{R}^n$, for each point of contact $p_i \in \mathcal{P}$, let
$\mathbf{k}_i$ be an outward pointing normal to $O$ at $p_i$ and $\pi_i$ the tangent
hyperplane $\mathbf{k}_i \cdot (x - p_i) = 0;$ then $O$ lies in the intersection of the half spaces $\mathcal{H}_i: \mathbf{k}_i \cdot (x - p_i) \leq 0.$\\
For immobilization  of $O$ it is necessary that $\mathcal{P}$ contains at least $n + 1$ points and that the intersection of all half-spaces $\mathcal{H}_i,$ is a bounded polytope, completely enclosing $O$; we assume this throughout below. Suppose that $\mathcal{P}$ contains precisely $n$ points and the corresponding normals are $\mathbf{k}_1,$ \ldots, $\mathbf{k}_n$. A displacement $\underline{u}$ of $O$ away from $p_i$ satisfies $\mathbf{k}_i \cdot \underline{u} < 0 .$ If the $n$ normals are linearly independent, then the system $\mathbf{k}_i \cdot \underline{u} = b_i, i = 1, \ldots, n$  has a unique solution for all choices of the $\{b_i\},$ in particular for all $b_i < 0.$ The displacement $\underline{u}$ moves $O$ away from all the $p_i$, so no penetration occurs. If the $\{\mathbf{k}_i\}$ are linearly dependent, then they span $W,$ a strict subspace of $\Bbl{R}^n$ and we simply choose $\underline{u}$ orthogonal to $W.$ Then $\mathbf{k}_i \cdot \underline{u} = 0$ for all $i,$ so $O$ 'slides' along each hyperplane $\pi_i$ and is again not immobilized by $\mathcal{P}$. Finally if $\mathcal{P}$ has less than $n$ points, $O$ is even less constrained than when $\mathcal{P}$ has $n$ points.

\noindent If the polytope  $\cap \mathcal{H}_i$ is unbounded, there is at least one direction $\underline{u}$ such that if any of the points  of $\mathcal{P}$ is displaced by an arbitrary positive multiple of $\underline{u},$ it remains within the polytope; then if $O$ is translated in this direction, it does not cut through any $\pi_i$, and so no penetration occurs.

\medskip 

\noindent Observe that if $O$ is a smooth convex body and the set $\mathcal{P}$ contains precisely $n + 1$ points $p_0, \ldots,  p_n$ the bounded polytope is a simplex $\Delta$. Thus a necessary condition for immobilization of $O$ is immobilization of the bounding simplex $\Delta;$ the converse is false, since owing to the curvature of the boundary of $O$ into the interior of $\Delta,$ a displacement of $\mathcal{P}$ causing at least one or more $p_i$ to penetrate $\Delta$ may not cause any of these $p_i$ to penetrate $O.$ In view of this necessity, understanding the immobilization of simplices is important for the study of immobilization of smooth convex bodies in general.

\medskip
\noindent Let $\mathcal{K}$ denote the set $\{\mathbf{k}_0, \ldots, \mathbf{k}_n\}$ of outward normals to $O$ at the points $p_i$. In order to achieve the bounded simplex $\Delta$, it is necessary and sufficient that 

\medskip
\be
\item each subset $\mathcal{K}-\{ \mathbf{k}_i \}$ of $n$ outward normals is linearly independent; 
\item in the unique (up to an overall scalar multiple) dependency relation \\ \mbox{$\sum_{i=0}^{n} \lambda_i \mathbf{k}_i = \mathbf{0}$} among the outward normals $\mathbf{k}_i$, all the coefficients $\lambda_i$ are non-zero  and of the same sign.
\ee

\medskip
\noindent These conditions ensure that there is no direction $\underline{u}$ in $\Bbl{R}^n$ such that $\underline{u} \cdot \mathbf{k}_i \leq 0$ for $i = 0, 1, \ldots, n;$ then there is no translation $\underline{u}$ of $O$ which does not cause $O$ to cut into at least one of the $\pi_i.$ By rescaling, so that each outward normal $\mathbf{k}_i$ is replaced by the outward normal $|\lambda_i| \mathbf{k}_i,$ the conditions 1 and 2 above may equivalently be replaced by
\be \item $\mathcal{K}-\{ \mathbf{k}_0 \}$  is linearly independent
\item $\sum_{i=0}^n \mathbf{k}_i = \mathbf{0}.$
\ee
It then follows that apart from an overall positive scaling factor (positive to ensure that normals remain outward pointing), $\sum_{i=0}^n \mathbf{k}_i = \mathbf{0}$ is the unique dependency among the outward normals  $\mathbf{k}_i$, and that apart from this overall scaling factor, all the $\mathbf{k}_i$ are now fixed. 

\noi We assume  1 and 2 below and in the next section construct a particular set of outward normals to $\Delta$ for which 1 and 2 are satisfied, and furthermore, the overall scaling factor is fixed. To do this, in \S 3, we initially change our view point: we start with the simplex $\Delta,$ defined by its vertices, and in Proposition \ref{eqn_2_1} derive a neat algebraic replationship between the set of vertices and a set of outward normals which satisfy 1 and 2 and which have a particular overall scaling. Then, for consistency with the approach taken in \S 1, in Proposition \ref{matrixK}, we show how starting with the convex body $O$ and the contact set $\mathcal{P}$ with its associated outward normals, we can construct $\Delta$ and rescale the normals to obtain the same relationship between vertices and normals as was found in Proposition \ref{eqn_2_1}.

\section{Matrix Description of the Simplex $\Delta$ and Contact set $\mathcal{P}$} \label{matrix}

\noi Let  $\Delta$ be a simplex in $\Bbl{R}^n$ having vertex set $\mathcal{V} = \{ v_0, v_1, \ldots, v_n \}$. Suppose the vertices are oriented so that $\vol(\Delta)$, the $n-$ dimensional volume of $\Delta,$ is positive. Let $[a_1, \ldots, a_n]$ denote the $n \times n$ matrix with columns $a_1, \ldots, a_n$. Then $\vol(\Delta) =  \frac{1}{n!}\det[v_1 - v_0, \ldots, v_n - v_0].$ 
To maintain symmetry, however, we regard each vertex $u$ with coordinates $(u_1, \ldots, u_n)$ as the point $\bar{u} = (1, u_1, \ldots, u_n) = (1, u)$ in the hyperplane $x_0 = 1$ in $\Bbl{R}^{n + 1}.$ Then letting $V$ be the $(n + 1) \times (n + 1)$ matrix $[\bar{v}_0, \bar{v}_1, \ldots, \bar{v}_n]$ it follows that $\vol(\Delta) = \frac{1}{n!}\det(V).$

\noi Let the face $F_i, i = 0, 1, \ldots, n$ of $\Delta$ be the $(n-1)$ dimensional simplex with vertex set $\mathcal{V}-\{v_i\}$, given a positive orientation. Let $\vol(F_i)$ denote the $(n-1)$ dimensional volume of $F_i$. Then if $h_i$ is the length of the altitude dropped from vertex $v_i$ to face $F_i$, regarding $\Delta$ as a cone with base $F_i,$ its volume is given by

\begin{equation}
\vol(\Delta) = \int_0^{h_i} \left(\frac{x}{h_i}\right)^{n-1} \vol(F_i)~dx = \frac{1}{n} h_i \vol(F_i). \label{eqn_2_1}
\end{equation}

\begin{prop} Let $\Delta$ be a simplex in $\Bbl{R}^n$ having  vertex set $\mathcal{V} = \{v_0, \ldots, v_n\}$ with an orientation such that the $n-$dimensional volume  of $\Delta$ is positive. There exists a set of outward pointing normals $\mathbf{k}_i$ to the faces $F_i$ of $\Delta$ satisfying:
\be 
\item $\sum_{i=0}^n \mathbf{k}_i = \mathbf{0},$ \item $|\mathbf{k}_i| = \vol(F_i)$, $i = 0, 1, \ldots, n.$\ee  \end{prop}

\noi {\bf Proof. } 
\noi Since $\vol(\Delta) > 0$, the matrix  $V$ is invertible, so there exists the $(n+1) \times (n+1)$ matrix $K = [\bar{\mathbf{k}}_0,  \bar{\mathbf{k}}_1, \ldots, \bar{\mathbf{k}}_n]$ such that 

\begin{equation}K^T V = V K^T = -n \vol(\Delta) I. \label{eqn_2_2} \end{equation}

\noi Analogous to the decomposition $\bar{v} = (1, v)$, we  write each $\bar{\mathbf{k}}_i = (\kappa_i, \mathbf{k}_i)$; the $n-$vectors $\bar{\mathbf{k}}_i$, $i = 0, 1, \ldots, n$ are the vectors we seek as we now show. Then from $K^TV = -n \vol(\Delta) I$ there follows
\begin{equation} \kappa_i + \mathbf{k}_i \cdot v_i =  - n \vol(\Delta),~~ 0 \leq i \leq n, \label{eqn_2_3} \end{equation}
\begin{equation} \kappa_i + \mathbf{k}_i \cdot v_j  = 0, ~~~0 \leq j \neq i \leq n, \label{eqn_2_4} \end{equation}

\noi while from the first row of $VK^T = -n \vol(\Delta) I,$ there follows

\begin{equation} \sum_{i=0}^n \kappa_i = -n \vol(\Delta), \label{eqn_2_5} \end{equation}
\begin{equation} \sum_{i=0}^n \mathbf{k}_i = \mathbf{0}. \label{eqn_2_6}\end{equation}

\noi From (\ref{eqn_2_4}) for each $i$ and each $j, l \neq i$

\begin{equation} \mathbf{k}_i \cdot (v_j - v_l) = 0.  \end{equation}


\noi Thus $\mathbf{k}_i$ is perpendicular to each edge of $F_i$ and is thus normal to $F_i.$

\medskip
\noi For any point $p_i$ of $F_i,$ since $v_j \in F_i,$ for all $j \neq i$, it follows that $p_i-v_j$ is parallel  to $F_i$ and hence  
\begin{equation}
\label{new_eqq1} \mathbf{k}_i \cdot (p_i - v_j) = 0
\end{equation} 
which, using (\ref{eqn_2_4}), gives for each $0 \leq i \leq n$ and for all $j \neq i$
\begin{equation}
\label{new_eqq2} \kappa_i = - \mathbf{k}_i \cdot p_i = - \mathbf{k}_i \cdot v_j .
\end{equation} 
From (\ref{eqn_2_3}) and (\ref{eqn_2_4}), for each $i$ and each $j \neq  i$, 

\begin{equation} \mathbf{k}_i \cdot (v_j - v_i) = n \vol(\Delta). \label{eqn_2_8} \end{equation}

\noi Now since $\mathbf{k}_i $ is perpendicular to $F_i$, the projection of $(v_j - v_i)$ along $\mathbf{k}_i$ has the length $h_i$ of the altitude from vertex $v_i$ to the face $F_i$. Thus (\ref{eqn_2_8}) implies $|\mathbf{k}_i| h_i =  n \vol(\Delta)$ and comparison with (\ref{eqn_2_1}) shows that $|\mathbf{k}_i| = \vol(F_i).$

\noi Finally by (\ref{eqn_2_8}) since $n \vol(\Delta) > 0, \mathbf{k}_i$ is outward pointing.

\medskip
\noi In this section, we started with $\Delta,$ specified by $\mathcal{V};$ this led to the faces $F_i$ and then,  via equation (\ref{eqn_2_2}), to the particular outward normals $\mathbf{k}_i$ of magnitude $\vol(F_i)$ encoded by the matrix $K$ in which the first row comprised elements $\kappa_i = -\mathbf{k}_i \cdot p_i$ with $p_i$ any point in $F_i$. We now show that we reach the same endpoint if we start, as in \S~\ref{westart}, with the points $p_i$ and the outward normals $\mathbf{k}_i$ which satisfy $\mathbf{k}_0 + \mathbf{k}_1 + \cdots + \mathbf{k}_n = \mathbf{0}$, $\mathbf{k}_1, \ldots, \mathbf{k}_n$ linearly independent and with a normalisation to be specified below.

\begin{prop} Let $O$ be a convex body in $\Bbl{R}^n$ and $\mathcal{P} = \{p_0, p_1, \ldots , p_n \}$ be a contact set of $O$. Let $\mathbf{k}_i$ be the outward pointing normals at $p_i$  and let these normals satisfy conditions 1 and 2 of \S 2, namely $\mathbf{k}_1, \ldots \mathbf{k}_n$ are linearly independent and $\sum_{i=0}^n \mathbf{k}_i = \mathbf{0}$. Then if for $i = 0, \ldots, n$, $\pi_i$ denotes the hyperplane $\mathbf{k}_i \cdot x =   \mathbf{k}_i \cdot p_i,$ there is a unique $n-$simplex $\Delta$ with vertices $v_j, j = 0, \ldots, n$ where $v_j$ is the intersection point of all the $\pi_i$ with $i \neq j;$ by orienting the set $\mathcal{P}$ suitably we can ensure that $\vol(\Delta) > 0$.
Furthermore, if we write $\kappa_j = - \mathbf{k}_j \cdot p_j$ and denote, as before,  $V = [\bar{v}_0, \ldots, \bar{v}_n]$, $K = [\bar{\mathbf{k}}_0,  \bar{\mathbf{k}}_1, \ldots, \bar{\mathbf{k}}_n]$ where $\bar{v}_j = (1, v_j), \bar{\mathbf{k}}_j = (\kappa_j, \mathbf{k}_j)$, then there is an overall positive rescaling of the outward normals $\{ \mathbf{k}_i \}$ so that $K$ and $V$ satisfy $K^TV = VK^T = -n \vol(\Delta) I,$ from which follow all the relations between normals and vertices established in Proposition \ref{eqn_2_1}. \label{matrixK} \end{prop}

\noi {\bf Proof. }
The vertex $v_j$ is the intersection point  of the $n$ hyperplanes $\pi_i:$\hspace{1cm}$ \; \mathbf{k}_i \cdot x =  \mathbf{k}_i \cdot p_i,$ $~~~0 \leq i \neq j \leq n.$ Since the $n$ normals $\mathbf{k}_i, i \neq j$ are linearly independent this defines $v_j$ uniquely and hence the simplex $\Delta$ with vertex set $\mathcal{V} = \{v_0, \ldots, v_n\}.$ Note that $\Delta$ and $\mathcal{V}$ are independent of the overall scaling of the normals since $\mathbf{k}_i$ appears linearly in both sides of the equation for $\pi_i.$ Thus for $j = 0, 1, \ldots, n, v_j$ satisfies
\begin{equation}
\mathbf{k}_i \cdot v_j = \mathbf{k}_i \cdot p_i, ~~~0 \leq i \neq j \leq n \label{eqn3p11} 
\end{equation}
and writing $\kappa_i = -\mathbf{k}_i \cdot p_i,$ there follows
\begin{equation}
\kappa_i + \mathbf{k}_i \cdot v_j = 0, ~~~~0 \leq i \neq j \leq n. \label{eqn3p12}
\end{equation}
Now summing (\ref{eqn3p12}) over all $i \neq j$ and using $\mathbf{k}_j = - \sum_{i\neq j} \mathbf{k}_i,$ we have
\begin{equation}
\kappa_i + \mathbf{k}_i \cdot v_i = \sum_{r=0}^n \kappa_r;~~~~~~ 0 \leq i \leq n \label{eqn3p13}
\end{equation}
and we note that (\ref{eqn3p12}), (\ref{eqn3p13}) also hold independently of the overall scaling of the normals.

\noi Now $\sum_{r=0}^n \kappa_r < 0:$ to see this, let $X$ be strictly inside $O.$ Then since $O$ is convex and the outward normal at each contact point $p_r$ is $\mathbf{k}_r$, we must have $(p_r - X) \cdot \mathbf{k}_r > 0, r = 0 , \ldots, n.$ Hence $\sum_{r=0}^n \kappa_r = - \sum_{r=0}^n \mathbf{k}_r \cdot p_r = - \sum_{r=0}^n \mathbf{k}_r \cdot (p_r - X) < 0.$
Defining vectors $\bar{v}_j, \bar{\mathbf{k}}_l$ and matrices $V, K$ as in the statement of the Proposition, equations (\ref{eqn3p12}) and (\ref{eqn3p13}) then give
$K^TV = \left( \sum_{r = 0}^n \kappa_r \right) I$. Since $\sum_{r=0}^n \kappa_r < 0,$ $K$ and $V$ are non-singular, so that $\vol(\Delta) = \frac{1}{n!}\det V \neq 0.$

\noi Now by re-orienting the $p_i$ if necessary (and thus, by (\ref{eqn3p11}), re-orienting the $v_i$) we can take $\vol(\Delta) > 0.$ Noting that $\sum_{r=0}^n \kappa_r < 0,$ we now do an overall positive scaling of the $\mathbf{k}_r$ (and hence all the $\kappa_r = - \mathbf{k}_r \cdot p_r$) so that 
$$\sum_{r=0}^n \kappa_r = - n \vol(\Delta).$$
We have thus now recovered (\ref{eqn_2_2}) 
$$K^TV = VK^T = -n \vol(\Delta) I.$$
Hence with this orientation of contact points and overall scaling of the outward normals, all the relations between normals and vertices established in Proposition \ref{eqn_2_1} follow.

\begin{rem} In practice to use (\ref{matrixK}) to obtain $V$ from $K,$
it is first necessary to compute $-n \vol(\Delta) = - n \det{V}/n!$ in terms of $K$ rather than $V.$ Taking determinants in (2) leads to $-n \vol(\Delta) = \left(- (n - 1)! \det K \right)^{\frac{1}{n}}$ and thus \\
$V = \left( - (n -1)! \det K \right)^{\frac{1}{n}} \left(K^T\right)^{-1},$ so giving the vertices. \end{rem}
 
\medskip
\noi We observe that 

\smallskip
\noi 1) If all the faces are projected orthogonally onto a single face $F_j$,  the sum of the projected volumes is zero (as is easily visualised in the cases $n = 2, 3$) which gives $\sum_{i = 0}^n \mathbf{k}_i \cdot \mathbf{k}_j = 0$. Thus by projecting onto faces $F_1, F_2, \ldots, F_n$ and using the linear independence of the corresponding normals, there follows $\sum_{i = 0}^n \mathbf{k}_i = \mathbf{0},$ thereby giving a geometric interpretation of (\ref{eqn_2_6}).
 
\medskip

\noi 2) From (\ref{eqn_2_4}) and (\ref{eqn_2_8}), for $j \neq i$, 
$$\kappa_i = -\mathbf{k}_i \cdot v_j = - \mathbf{k}_i \cdot (v_j - 0)$$
is $-n \vol(\Delta_i)$ where $\Delta_i$ is the simplex with vertices of $F_i$ together with the origin. As the origin is not in general a point of $\Delta$, the quantities $\kappa_i$ do not play a significant role later.

\medskip
\noi The contact points $p_i$ lie in the interior of each face $F_i$. Thus each $p_i$ is a linear combination of the vertices of $F_i$ with coefficients that are non-negative and sum to unity. Define $\bar{p}_i = (1, p_i)$ and $P = [\bar{p}_0, \bar{p}_1, \ldots, \bar{p}_n ]$;  we then have 
$$P = V \Lambda$$
\noi where $\Lambda$ is an $(n+1) \times (n+1)$ matrix whose columns list the coefficients of $\bar{p}_j$ in terms of $\bar{v}_0, \bar{v}_1, \ldots, \bar{v}_n.$ The matrix $\Lambda = (\lambda_{ij})$ then has the following properties: $\lambda_{ii} = 0$ since $v_i \notin F_i;$ $\lambda_{ij} > 0$ for all $i \neq j$; $\sum_{i=0}^n \lambda_{ij} = 1$. $\Lambda$ is thus a stochastic matrix with the additional property that the diagonal entries are zero. These matrices enjoy certain properties  which we exploit later. The matrix $\Lambda$ thus provides an efficient encoding of the set of contact points.

\section{The Penetration Function} \label{penetration}
The conditions for immobilization of  $\Delta$ are now recast in terms of an extremal problem. Let $\Bbl{E}_n$ denote the group of all rigid motions of Euclidean  space $\Bbl{R}^n$. Given the simplex described by the matrix $V$ (and thus its normals $K$) and the contact points by the matrix  $P$ define the penetration function $\Phi : \Bbl{E}_n \ra \Bbl{R}$ by

\begin{equation}\Phi (g) = \sum_{i=0}^n (g(p_i) - p_i) \cdot \mathbf{k}_i. \label{eqn_3_1} \end{equation}

\noi Then $\Phi$ varies continuously with $g$ and $-\Phi(g)$ measures the total amount of normal penetration into $\Delta$ by the points $p_i \in \mathcal{P}$ under the action of $g$ in $\Bbl{E}_n$ (weighted by the volumes of the faces).

\medskip
\noi Now each $g \in \Bbl{E}_n$ has a unique decomposition $g = tr$ where $t \in T_n,$ the group of translations of $\Bbl{R}^n$ and $r \in SO(n),$ the group of orientation preserving rotations of $\Bbl{R}^n$ about the origin, \cite{ONEILL}.
 For each $t \in T_n$ and $r \in SO(n)$, $r^{-1}tr$ is also a translation from which it follows that $T_n r = r T_n$ for all $r \in SO(n),$ so that $T_n$ is a normal subgroup of $E_n$ and $SO(n)$ is the quotient group.
 
\noi An element $g = tr \in \Bbl{E}_n$ has a convenient matrix representation. Let $t_{\ul{a}}$ denote the translation $x \mapsto x + \ul{a}$ and let $r$ have matrix $R.$ Then $g(x) = R x + \ul{a}$. As in \S~\ref{matrix}, let $\bar{x} = (1, x)$ and $\bar{a} = (1, \ul{a}).$ Then the $(n+1) \times (n+1)$ matrix $G$ with top row $(1, 0, 0, \ldots, 0)$ and lower $n \times (n+1)$ array $[\underline{a}, R]$ satisfies

 \begin{equation} G\bar{x}= \left[ \begin{array}{cc}\displaystyle{1} & \displaystyle{\underline{0}^T} \\ \displaystyle{\underline{a}} & \displaystyle{R} \end{array}\right] ~ \left[  \begin{array}{c} 1\\ x  
\end{array}  \right] =  \left[ \begin{array}{c} 1\\ gx  
\end{array} \right] = \ol{g(x)}. \label{eqn_3_2} \end{equation}

\noi With this formulation of $g = t_{\ul{a}} r$, it follows that if $g$ varies continuously so do $t_{\ul{a}}$ and $r$, and conversely. Thus the map from $\Bbl{E}_n$ to $T_n \times SO(n)$ given by $g = tr \mapsto (t, r)$ and its inverse are continuous.

\begin{lem} \label{lemma1}
\noi $\Phi$ is well defined on $\Bbl{E}_n/T_n = SO(n)$. \end{lem}

\noi {\bf Proof. }
Let $r \in SO(n);$ we show that $\Phi$ is constant on $T_nr$. For $g = t_{\ul{a}}r$
$$\Phi(t_{\ul{a}} r) = \sum_{i=0}^n \left( t_{\ul{a}}r p_i - p_i \right) \cdot \mathbf{k}_i = \sum_{i=0}^n(rp_i + \underline{a} - p_i) \cdot \mathbf{k}_i = \Phi(r)$$
\noi since $\sum_{i=0}^n \mathbf{k}_i = \mathbf{0}.$

\noi We thus  regard $\Phi$ as a continuous map from $SO(n)$ to $\Bbl{R}.$

\medskip
\begin{lem} \label{lemma2}
For each $r \in SO(n)$, there is a unique $t(r) \in T_n$ and hence $g(r) = t(r) r \in \Bbl{E}_n$ so that $\Phi(g(r)) = \Phi(r)$ and 
\begin{equation} \left(g(r) p_i - p_i\right) \cdot \mathbf{k}_i = \frac{\Phi(r)}{n + 1}, ~ i = 0, 1, \ldots, n. \label{eqn_3_3} \end{equation} In particular $g(r) = I \in \Bbl{E}_n$ if and only if $r = I \in SO(n).$ Furthermore $g(r)$ varies continuously with $r.$ \end{lem}

\begin{rem} Lemma \ref{lemma2} states that for each rotation in $SO(n)$ (and hence from the decomposition $g = t r$, for every rigid motion), there is a unique translation $t(r)$ so that the normal penetrations into $\Delta$ at the $p_i$ caused by $t(r) r$ are all equal. \end{rem}

\noi {\bf Proof. }
Equation (\ref{eqn_3_3}) holds if and only if there is a unique $\underline{a}$ so that 
\begin{equation}
(r p_i + \underline{a} - p_i) \cdot \mathbf{k}_i = \frac{\Phi(r)}{n + 1}, ~ i = 0, 1, \ldots, n. \label{eqn_3_4}
\end{equation}
\noi Since $\{ \mathbf{k}_1, \mathbf{k}_2, \ldots, \mathbf{k}_n \}$ are linearly independent, $\underline{a}$ is uniquely specified by the $n$ equations
\begin{equation}
(r p_i + \underline{a} - p_i) \cdot \mathbf{k}_i = \frac{\Phi(r)}{n + 1}, ~ i = 1, \ldots, n \label{eqn_3_5}
\end{equation}
\noi But then from Lemma \ref{lemma1}
\begin{eqnarray*} (r p_0 + \underline{a} - p_0) \cdot \mathbf{k}_0 & =  \sum_{i=0}^n (r p_i + \underline{a} - p_i) \cdot \mathbf{k}_i - \sum_{i=1}^n (r p_i + \underline{a} - p_i) \cdot \mathbf{k}_i \\
& =  \Phi(r) - n \frac{\Phi(r)}{n+1} = \frac{\Phi(r)}{n+1} 
\end{eqnarray*}
\noi which gives (\ref{eqn_3_4}). \\
When $r = I \in SO(n)$, then by (\ref{eqn_3_1}), $\Phi(I) = 0$ so that $\underline{a} = \underline{0}$ by (\ref{eqn_3_4}) and hence $g(r) = I \in \Bbl{E}_n$. Conversely if $r \neq I \in SO(n),$ then by (\ref{eqn_3_2}), $t_{\ul{a}} r \neq I \in \Bbl{E}_n$ for any $\underline{a}$, so in particular $g(r) \neq I.$ The system (\ref{eqn_3_4}) for $\underline{a}$ can be written

\begin{equation}
\mathbf{k}_i \cdot \underline{a} = b_i, ~ i = 0, 1, \ldots, n, \label{eqn_3_6}
\end{equation}

\noi where $b_i = \frac{\Phi(r)}{n+1} - (r p_i - p_i) \cdot \mathbf{k}_i.$ Using the notation of \S~\ref{matrix} and writing $\tilde{\underline{a}} = (0, \underline{a})$ and $\overline{b} = (b_0, b_1, \ldots, b_n),$ the system (\ref{eqn_3_6}) is 
$$K^T \tilde{\underline{a}} = \overline{b}$$ which has solution
$$\tilde{\underline{a}} = - \frac{V \ol{b}} {n \vol(\Delta)} $$ from which, invoking the continuity of $\Phi(r)$ with $r$, it follows that $\underline{a}$ and hence $t(r)$ and $g(r)$
 all depend  continuously on $r \in SO(n).$

\begin{lem} \label{lemma3}
For each $t \in T_n, \Phi(t) = 0.$ If $t \neq I \in \Bbl{E}_n$, then under the action of $t$ at least one point $p_i$ penetrates $\Delta.$ \end{lem}

\noi {\bf Proof. }
Since $\sum_{i=0}^n \mathbf{k}_i = \mathbf{0},$ $$\Phi(t_{\ul{a}}) = \sum_{i=0}^n (p_i + \underline{a} - p_i) \cdot \mathbf{k}_i = \underline{a} \cdot \sum_{i=0}^n \mathbf{k}_i = 0.$$
If $t_{\ul{a}} \neq I$, then $\underline{a} \neq 0$. Since $\{\mathbf{k}_1, \ldots,\mathbf{k}_n \}$ are linearly independent at least one $\underline{a} \cdot \mathbf{k}_i, ~i = 1, \ldots, n$ is non-zero and since $\sum_{i=0}^n \mathbf{k}_i = \mathbf{0}$  at least one $\underline{a} \cdot \mathbf{k}_i, ~i = 0, \ldots, n$ is strictly negative. Hence $(t_{\ul{a}} p_i - p_i) \cdot \mathbf{k}_i = \underline{a} \cdot \mathbf{k}_i < 0$ for at least one $i,$ so that $p_i$ penetrates $\Delta$ under the action of $t_{\ul{a}}.$ 

\medskip
\begin{prop} \label{prop4}
The set $\mathcal{P}$ immobilizes $\Delta$ if and only if $\Phi: SO(n) \ra \Bbl{R}$ has a strict local maximum at $I.$ \end{prop}

\noi {\bf Proof. }
Let $\mathcal{P}$ immobilize $\Delta.$ Then each rigid motion $g \neq I$ in a sufficiently small neighbourhood $N$ of $I$ in $\Bbl{E}_n$ causes at least one $p_i$ to penetrate $\Delta.$ By the continuity of $r \mapsto g(r)$ there is a neighbourhood $N'$ of $I \in SO(n)$ so that $r \in N'$ implies $g(r) \in N.$  Now suppose that $\Phi$ does not have a strict local maximum at $I \in SO(n).$ Then each neighbourhood of $I \in SO(n)$, and in particular $N'$, 
contains a rotation $r \neq I$ such that $\Phi(r) \geq 0$. For this $r$, $g(r) \in N$, and by Lemma \ref{lemma2}, $g(r) \neq I$ and $(g(r) p_i - p_i ) \cdot \mathbf{k}_i = \frac{\Phi(r)}{n+1} \geq 0$ for $i = 0, 1, \ldots, n,$ so that $g(r)$ causes no $p_i$ to penetrate $\Delta;$ this gives the required contradiction.

\medskip
\noi Conversely let $\Phi$ have a strict local maximum at $I \in SO(n).$ Then $\Phi(r) < 0$ for all $r \neq I$ in a sufficiently small neighbourhood $N'$ of $I$ in $SO(n).$ Now using the continuity of $tr \mapsto (t, r),$ consider a rigid motion $g = t r \neq I$ in a sufficiently small neighbourhood of $I$ in $\Bbl{E}_n$ so that $r \in N'.$ If $r = I,$ then $g = t \neq I$ and hence by Lemma \ref{lemma3}, at least one $p_i$ penetrates $\Delta;$ if $r \neq I,$ then $\Phi(g) = \Phi(r) < 0$ so that at least one term $(g(p_i) - p_i) \cdot \mathbf{k}_i < 0$ so that $p_i$ penetrates $\Delta.$ Thus $\mathcal{P}$ immobilizes $\Delta.$ 

\section{Conditions for a Maximum of the Penetration Function} \label{conditions}
Given $\Delta$ and $\mathcal{P}$, we now examine conditions under which $\Phi$ has a strict local maximum at $I \in SO(n).$ Regarding $g \in SO(n)$ as an $n \times n$ matrix and noting that $p_i \cdot \mathbf{k}_i$ can be identified as the trace of the $n \times n$ matrix $\mathbf{k}_i p_i^T,$ it follows that
\begin{equation}
\sum_{i=0}^n g(p_i) \cdot \mathbf{k}_i = \sum_{i=0}^n \tr(\mathbf{k}_ip_i^T g^T) = \sum_{i=0}^n \tr(g^T \mathbf{k}_i p_i^T) = \tr(g^T A) = \tr(A^T g), \label{eqn_4_1}
\end{equation}
\noi where $A$ is the matrix
\begin{equation}
A = \sum_{i=0}^n \mathbf{k}_i p_i^T \label{eqn_4_2}
\end{equation}
\noi which depends solely on the geometry of $\Delta$ and the set of contact points $\mathcal{P}.$ Hence
\begin{equation}
\Phi(g) = \tr(A^T(g - I)) \label{eqn_4_3}
\end{equation}
and we now consider conditions on $A$ so that $\Phi$ has a strict local maximum at $I.$ 

\medskip
\noi Now each $g \in SO(n)$ may be written as $$g = \exp S = \sum_{k \geq 0} \frac{S^k}{k !}$$
where $S$ is a skew symmetric matrix, and where $g = I$ corresponds to $S = 0.$ Hence $\Phi(g) = \Psi(S)$ where 
\begin{equation}
\Psi(S) = \tr(A^TS) + \frac{1}{2!} \tr(A^T S^2) + \cdots   \label{eqn_4_4}
\end{equation}
\noi We regard (\ref{eqn_4_4}) as a power series about $S = 0$ and examine it for a strict local maximum at $S = 0.$
By definition $\Psi$ has a strict local maximum at $S = 0$ if $\Psi(S) < 0$ for all $S \neq 0$ in some neighbourhood of $S = 0.$ In the following we make repeated use of the result that there exists a sufficiently small neighbourhood $N$ of $S = 0$ such that for all $S \neq 0$ in $N,$ the sign of $\Psi(S)$ is the same as the sign of the first non-zero term in the power series (\ref{eqn_4_4}).

\medskip
\begin{lem} \label{lemma1,2}
For $\Psi$ to have a local maximum at $S = 0,$ it is necessary that for all skew symmetric matrices $S$, $\tr(A^T S) = 0.$ \end{lem}

\noi {\bf Proof. } This is just the usual condition to make $S = 0$ a stationary point of $\Psi$, which is necessary for an extreme value at $S = 0$. 

\medskip
\noi We now examine further conditions to ensure sufficiency of a strict local maximum at $S = 0.$

\begin{lem} \label{lemma2,2}
For $\Psi$ to have a strict local maximum at $S = 0,$ it is sufficient that for all skew-symmetric matrices $S \neq 0$: \\
(i) $\tr(A^TS) = 0$; \\
(ii) $\tr(A^TS^2) < 0.$ \end{lem}

\noi {\bf Proof. }
Assuming  (i) and (ii), it follows from (i) that
\begin{equation}
\Psi(S) = \frac{1}{2!}\tr(A^TS^2) + \frac{1}{3!}\tr(A^TS^3) + \cdots \label{eqn_4_5}
\end{equation}

\noi For all  $S \neq 0$ in a neighbourhood $N,$ the sign of $\Psi(S)$ is the same as that of $\tr(A^TS^2).$ Since by (ii), $\tr(A^TS^2) < 0$ for all $S \neq 0$ in $N$, it follows that $\Psi(S) < 0$ for all $S \neq 0$ in $N$ and hence $\Psi$ has a strict local maximum at $S = 0.$

\medskip
\noi We now examine the implications for $A$ to ensure that conditions (i) and (ii) hold.

\begin{lem} \label{lemma3,2}
$\tr(A^TS) = 0$ for all skew-symmetric matrices $S$ if and only if $A$ is symmetric. \end{lem}

\noi {\bf Proof. }
Let $S^{(ij)}$ be the skew symmetric matrix with a 1 in the position $(i,j), -1$ in position $(j, i) $ and 0 everywhere else. Then $\tr\left(A^T S^{(ij)}\right) = a_{ji} - a_{ij} = 0.$ The converse follows as each skew symmetric $S$ can be expressed as $S = \sum_{i < j}c_{ij}S^{(ij)}.$ 

\begin{defn} \label{dfn_apd}
A real symmetric matrix $A$ is almost positive definite if the sum of every pair of eigenvalues is positive. \end{defn}
\noi Here, an eigenvalue can only be added to itself if it is a repeated eigenvalue. This condition is weaker than positive definiteness.

\medskip
\begin{lem} \label{lemma4}
For a real symmetric matrix $A$, $\tr(A^TS^2) < 0$ for all skew symmetric matrices $S \neq 0$ if and only if $A$ is almost positive definite. \end{lem}

\noi {\bf Proof. }
We have $\tr(A^TS^2) = \tr(SAS) = -\tr(S^TAS).$ $A$ can be diagonalised as $P^TAP = D$ where $P$ is orthogonal. Hence
$$\tr(S^TAS) = \tr(S^TPDP^TS) = \tr(P^TS^TP \:D \: P^TSP).$$
\noi Now $S$ is skew-symmetric if and only if $P^TSP$ is skew-symmetric, hence \\ $\tr(S^TAS) > 0$ for all skew  $S \neq 0$ if and only if $\tr(S^TDS) > 0$ for all skew $S \neq 0.$ Choosing now $S = S^{(ij)},$ then $S^{(ij)}S^{(ij)T}$ is diagonal with a 1 in position $i$ and $j$ and 0 elsewhere. Hence 
$$\tr(S^{(ij)T}D\, S^{(ij)}) = \tr(S^{(ij)}S^{(ij)T}D) = \lambda_i  + \lambda_j$$ so that if $S = \sum_{i < j}c_{ij} S^{(ij)},$ then 
$$\tr(S^T D S) = \sum_{i < j} c_{ij}^2 (\lambda_i + \lambda_j) > 0$$ if and only if each pair $\lambda_i + \lambda_j > 0,$ which gives the result. 

\medskip

\noi Combining Lemmas \ref{lemma2,2}, \ref{lemma3,2} and \ref{lemma4} we thus have:

\medskip
\begin{cor} \label{cor5}
$\Psi$ has a strict local maximum at $S = 0$ if (i) $A$ is symmetric and (ii) $A$ is almost positive definite. \end{cor}

\medskip
\noi We want to examine the extent to which the conditions (i) and (ii) are necessary. Given that (i) holds, from (\ref{eqn_4_5}) considering what happens with a real function of a real variable, one could envisage $\Psi$ having a strict quartic maximum at $S = 0$, given by
$$\Psi(S) = \frac{1}{4!}\tr(A^TS^4) + \cdots$$ in which $\tr(A^TS^2) = 0$, $\tr(A^TS^3) = 0$ and $\tr(A^TS^4) < 0$. This however does not happen; it transpires that Corollary \ref{cor5} has a converse. Thus we have

\medskip
\begin{prop} \label{prop6}
$\Psi$ has a strict local maximum at $S = 0$ if and only if: \\
(i) $A$ is symmetric; and \\
(ii) $A$ is almost positive definite. \end{prop}

\noi {\bf Proof. }
The if part is given by Corollary \ref{cor5}. \\
Now (i) is necessary by Lemma \ref{lemma1,2} and Lemma \ref{lemma3,2}. To establish the necessity of (ii) we now assume that $\Psi$ has a strict local maximum at $S = 0$ and that $A$ is symmetric, in which case $\Psi(S)$ is given by (\ref{eqn_4_5}).

\noi We complete the proof by contradiction, so we assume that $A$ is not almost positive definite. Then by Lemma \ref{lemma4}, there is a skew $S_0 \neq 0$ so that $\tr(A^T S_0^2) \geq 0.$ We  now consider two cases: 

\noi (a) If $\tr(A^TS_0^2) > 0,$ then replacing $S_0$ by $\eps S_0$ where $\epsilon > 0$ is   sufficiently small that $\epsilon S_0 \in N,$ the neighbourhood of $S = 0$ defined after (\ref{eqn_4_4}), it follows that $\Psi(\epsilon S_0) > 0$, giving the required contradiction;

\noi (b) If there is no skew $S \neq 0$ such that $\tr(A^TS^2) > 0,$ then $\tr(A^TS_0^2) = 0.$ Using the same arguments as those in Lemma \ref{lemma4} we may take $A$ to be diagonal; also, since there is no $S \neq 0$ such that $\tr(A^TS^2) > 0,$ there is no pair of eigenvalues $\lambda_i$, $\lambda_j$ with $\lambda_i + \lambda_j < 0;$ hence $\lambda_i + \lambda_j \geq 0$ for all pairs of eigenvalues. Next if $S_0 = \sum_{i < j} c_{ij} S^{(ij)}$ then from $\tr(A^TS_0^2) = 0$ there follows $\sum_{i < j} c_{ij}^2(\lambda_i + \lambda_j) = 0$ so that $\lambda_i + \lambda_j = 0$ for all pairs $i, j$ for which $c_{ij} \neq 0.$ We can thus assume $\lambda_1 + \lambda_2 = 0.$ Then for $S = S^{(12)},$ we have that $S^2$ is diagonal with entries $-1, -1, 0, 0, \ldots, 0$ from which follows $\tr(A^TS^{2r}) = (-1)^r(\lambda_1 + \lambda_2) =  0$ for all $r \geq 1$ and since $S^{2r + 1}$ is skew $\tr(A^TS^{2r+1}) = 0$ for all $r \geq 0.$ Hence for $S = S^{(12)} \neq 0,$ we have $\Psi(S) = 0$, again contradicting the existence of a strict local maximum at $S = 0.$

\medskip
\noi Proposition \ref{prop6} gives necessary and sufficient conditions for the immobilization of $\Delta.$ These are couched as algebraic conditions on the matrix $A$: we recall that $A$ depends solely on the geometry of $\Delta$ and the contact set $\mathcal{P}$.\\
We now establish a relation between the matrix $A$ and the matrix $\Lambda,$ introduced in \S~\ref{matrix}, and which like $A$ also only depends on the geometry of $\Delta$ and $\mathcal{P}$; we then use the special properties of $\Lambda$ to show that symmetry of $A$ alone suffices for immobilization of $\Delta$ in the cases $n = 2, 3,$ but not for $n \geq 4.$

\medskip
\begin{lem} \label{lemma7}
The eigenvalues of the $(n+1) \times (n+1)$ matrix $-n \vol(\Delta) \Lambda$ are precisely those of the $n \times n$ matrix $A$ together with the value $-n \vol(\Delta).$ \end{lem}

\medskip
\noi {\bf Proof. }
The matrix $A = \sum_{i=0}^n \mathbf{k}_i p_i^T$ is the product of the $n \times (n+1)$ matrix $[\mathbf{k}_0, \mathbf{k}_1, \ldots, \mathbf{k}_n]$
with the transpose of the $n \times (n+1)$ matrix $[p_0, p_1, \ldots, p_n].$ Reintroducing the $(n+1) \times (n+1)$ matrices $K = [\bar{\mathbf{k}}_0,  \bar{\mathbf{k}}_1, \ldots \bar{\mathbf{k}}_n]$ and $P = [\bar{p}_0, \bar{p}_1, \ldots, \bar{p}_n]$ we have
\begin{equation}
KP^T = \left[\begin{array}{cccc} \kappa_0 & \kappa_1 & \cdots & \kappa_n \\ \mathbf{k}_0 & \mathbf{k}_1 & \cdots & \mathbf{k}_n \end{array} \right] \left[ \begin{array}{cc} 1 & p_0^T \\ 1 & p_1^T \\ \vdots & \vdots \\ 1 & p_n^T \end{array}\right] = \left[ \begin{array}{cc} -n \vol(\Delta) & \underline{b}^T \\ \mathbf{0} & A \end{array} \right], \label{eqn_4_6}
\end{equation}

\noi where we have used $\sum \kappa_j = -n \vol(\Delta), ~\sum \mathbf{k}_j = \mathbf{0}$ and where $\underline{b} = \sum_j \kappa_j p_j$ plays no further role in the discussion. Thus the eigenvalues of $KP^T$ are precisely those of $A$ together with a further eigenvalue $-n \vol(\Delta).$\\
But from \S~\ref{matrix}, $P = V \Lambda$ and $KV^T = -n \vol(\Delta) I$ so that 
$$KP^T = K \Lambda^T V^T = -n \vol(\Delta) K \Lambda^T K^{-1}.$$ Hence $KP^T$ is similar to $-n \vol(\Delta) \Lambda^T$ which gives the result. 

\bigskip
\noi Now from \S~\ref{matrix}, $\Lambda$ is a stochastic matrix with diagonal entries all zero and all off-diagonal entries positive. The following proposition lists the possibilities for the eigenvalues of such a matrix (the proof is actually for $\Lambda^T$  which has the same eigenvalues as $\Lambda$).

\medskip
\begin{prop} \label{prop8}
Let $\Lambda = (\lambda_{ij})$ be a $m \times m$ matrix such that:\\
(i) $\sum_{j=1}^m \lambda_{ij} = 1$ for $1 \leq i \leq m;$ \\
(ii) $\lambda_{ii} = 0$ for all $1 \leq i \leq m;$\\
(iii) $\lambda_{ij} > 0$ for all $1 \leq i \neq j \leq m$.\\
Then for all $m \geq 3:$
\begin{itemize}
\item 1 is an eigenvalue of $\Lambda$ with eigenvector $(1, 1, \cdots, 1)^T$ \item all other eigenvalues lie strictly within the unit circle.
\end{itemize} \end{prop}

\noi {\bf Proof. }
We consider vectors in $\Bbl{C}^m$ with the supremum norm $||\underline{z}|| = \max\{|z_k|: 1 \leq k \leq m\}.$ Then for $k = 1, 2, \ldots, m$
$$|(\Lambda \underline{z})_k| = \left| \sum_{j=1}^m \lambda_{kj} z_j \right| \leq \sum_{j=1}^m |\lambda_{kj}||z_j| \leq \sum_{j=1}^m \lambda_{kj} ||\underline{z}|| = ||\underline{z}||.$$ 
Hence if $\lambda$ is an eigenvalue of $\Lambda$ with eigenvector $\underline{z}$

$$|\lambda|||\underline{z}|| = ||\lambda \underline{z}|| = ||\Lambda \underline{z}|| = \max\{|(\Lambda \underline{z})_k|: 1 \leq k \leq m\}  \leq ||\ul{z}||,$$
\noi so that $|\lambda| \leq 1.$

\medskip
\noi If $\underline{z} = (1, 1, \ldots, 1)^T$, then $(\Lambda \underline{z})_k = \sum_{j=1}^m \lambda_{kj} = 1$ so that $\Lambda \underline{z} = (1, 1, \ldots, 1)^T$ and hence $\underline{z} = (1, 1, \ldots, 1)^T$ is an eigenvector corresponding to $\lambda = 1.$ 

\medskip
\noi Now let $\Lambda \underline{z} = \mu \underline{z}$ where $|\mu| = 1$ and let $|z_{k^*}|$
be maximal among $|z_k|, ~1 \leq k \leq m.$ Then

\begin{equation}
|z_{k^*}| = |\mu z_{k^*}| = |(\Lambda \underline{z})_{k^*}| = \left|\sum_{j=1}^m \lambda_{k^*j} z_j\right|  \leq \sum_{j=1}^m \lambda_{k^*j} | z_j|  \label{eqn_4_7}
\end{equation}

\noi and using $\sum_j \lambda_{k^*j} = 1,$ it follows that 
$$0 \leq \sum_{j=1}^m \lambda_{k^*j} \left(|z_j| - |z_{k^*}|\right).$$ 

\noi But since $\lambda_{k^*j} > 0$ and $|z_j| - |z_{k^*}| \leq 0$ for all $j \neq k^*,$ it follows $|z_j| = |z_{k^*}|$ for all $j \neq k^*.$ Thus the components $z_j$ of eigenvector  $\underline{z}$ all lie on a circle of radius $|z_{k^*}|$ in $\Bbl{C},$ and for $m \geq 3,$ unless all the $z_k$ are the same, the point $\sum_{j=1}^m \lambda_{k^*j} z_j$ lies strictly inside the convex hull of the $\{z_j\}$ and so has modulus strictly less than $|z_{k^*}|$ contradicting the equality $|z_{k^*}| = |\sum_{j=1}^m \lambda_{k^*j} z_j|$ in (\ref{eqn_4_7}) above. Hence the eigenvector $\underline{z}$ is  a multiple of $(1, 1, \ldots, 1)^T$ and thus corresponds to $\lambda = 1;$ thus all other eigenvalues $\lambda$ satisfy $|\lambda| < 1.$  

\medskip
\begin{rem}The final argument relies on the fact that for $m \geq 3, \Lambda$ has elements satisfying $0 < \lambda_{ij} < 1;$ this is not true when $m = 2$ in which case 
$$ \Lambda = \left[\begin{array}{cc} 0 & 1 \\ 1 & 0\end{array} \right]$$
 and then $\Lambda$ has an eigenvalue of $-1$ corresponding to eigenvector $(-1, 1)^T.$ \end{rem}

\medskip
\noi This proposition leads to a difference in the sufficiency conditions for immobilization between the cases $n \leq 3$ and $n > 3.$

\medskip
\begin{prop} \label{prop9}
When $n \leq 3$ a sufficient condition for immobilization is simply the symmetry of $A.$ \end{prop}

\noi {\bf Proof. }
We show that for $n \leq 3,$ the fact that the points $p_i$ are internal to the faces $F_i$ together with symmetry of $A$ implies that $A$ is almost positive definite. The result then follows from Corollary \ref{cor5}.

\noi If $A$ is symmetric its eigenvalues are real and following Lemma \ref{lemma7}, we write them as $-n \vol(\Delta)\lambda_1$, $\ldots$, $ -n \vol(\Delta)\lambda_n,$ where $\lambda_1, \ldots, \lambda_n$ are the eigenvalues of $\Lambda$ other than unity. Since $\Lambda$ is tracefree, $\lambda_1 + \lambda_2 + \cdots + \lambda_n + 1 = 0.$ For $n = 2,$ $\lambda_1 + \lambda_2 = -1$ so that $-n \vol(\Delta)\lambda_1 - n \vol(\Delta)\lambda_2 = n \vol(\Delta)$ and so $A$ is almost positive definite. For $n = 3,$ $\lambda_1 + \lambda_2 + \lambda_3 + 1 = 0,$ so, for example $\lambda_1 + \lambda_2 = -(1 + \lambda_3) < 0$ since $|\lambda_3| < 1$ by Proposition \ref{prop8} and $\lambda_3$ is  real, similarly for each pair of eigenvalues and hence $A$ is almost positive definite. 

\noi This method clearly fails for $n \geq 4$ and indeed the reduced sufficiency condition does not hold for $n \geq 4$ as the example in Appendix A shows.

\medskip
\noi For the case $n = 3$ Proposition \ref{prop9} was implied by the work of \cite{BFMM} and for the case $n = 2$, it was implied by the work of \cite{CSU}. In both cases different methods were used.  The use here of stochastic matrices which are widely studied seems advantageous.

\section{Geometric Interpretation of Conditions for \\ Immobilization} \label{geometric}
We finally attempt to give a geometric description of the conditions for immobilization of $\Delta$. The ideal result would be to give a complete geometric description of the elements of $\mathcal{S}$, the set of all contact sets $\mathcal{P}$ which immobilize $\Delta.$ This is possible when  $n = 2,$ (see \cite{CSU}) but already becomes awkward for $n = 3$. Instead, we settle for something less. We are able to give a good geometric description of some of the elements of $\mathcal{S}$ and one which we call $\mathcal{G}$ stands out as being in some sense optimal. We will also describe a geometric process, which starting with $\mathcal{P} = \mathcal{G}$ displaces the point $p_i,$ each within its face $F_i,$ to reach other elements of $\mathcal{S}$. In this process, symmetry of $A$ is maintained, but it is harder to keep track of almost positive definiteness.

\medskip

\noi Our procedure is as follows: we initially construct some special $n \times (n+1)$ matrices $[p_0, p_1, \ldots, p_n]$ of contact points which ensure that $A$ is symmetric, but where we temporarily relax the condition $p_i \in F_i.$ From the linearity of $A$ with respect to the $p_i$ in $A = \sum \mathbf{k}_i p_i^T,$ it follows that any linear combination of these special matrices also yields a matrix $A$ which is symmetric, and we choose the linear combinations to ensure that $p_i \in F_i;$ we will also examine almost positive definiteness of such matrices $A.$

\noi The following special sets $[p_0, p_1, \ldots, p_n]$ guarantee symmetry of $A:$

\medskip
\noi (i) $[t_0\mathbf{k}_0, t_1\mathbf{k}_1, \ldots, t_n\mathbf{k}_n]$ where $t_i, i = 0, 1, \ldots, n$ are arbitrary. Then $A = \sum_{i=0}^n t_i \mathbf{k}_i\mathbf{k}^T_i$ is clearly symmetric and positive definite if all $t_i  > 0.$\\
(ii) $[\underline{z}, \underline{z}, \ldots, \underline{z} ]$ where $\underline{z} \in \Bbl{R}^n$ is arbitrary. Then $A = \sum \mathbf{k}_i \underline{z}^T$ which is the zero matrix.\\
(iii) $[v_0, \ldots, v_n]$. This corresponds to the $(n+1) \times (n+1)$ matrix \\ $V = [\bar{v}_0, \bar{v}_1, \ldots, \bar{v}_n]$  for which $KV^T = -n \vol(\Delta) I_{n+1}$ and then by (\ref{eqn_4_6}), $A = -n \vol(\Delta) I_n.$\\
(iv) $[\mathbf{k}_1, \mathbf{k}_0, 0, \ldots, 0]$  gives $A = \mathbf{k}_0 \mathbf{k}_1^T + \mathbf{k}_1\mathbf{k}_0^T$ which is clearly symmetric: similarly any interchange of $\mathbf{k}_i$ and $\mathbf{k}_j,$ $i \neq j$ with all other entries 0. 

\medskip
\noi We now use these special sets to construct some immobilizing sets.

\subsection{Centroids of the Faces} \label{subsec:centroids}

\noi Let $g = \frac{1}{n} \sum_{i=0}^n v_i.$ Then $g_i = g - \frac{v_i}{n}$ is the centroid of the face $F_i$. Let $G = [\bar{g}_0, \bar{g}_1, \ldots, \bar{g}_n] $ be the $(n+1) \times (n+1)$ matrix of the contact set of centroids. Then for this set

$$KG^T = K [\bar{g}, \bar{g}, \ldots, \bar{g}]^T - \frac{1}{n} KV^T = 0 + \vol(\Delta)I_{n+1}$$

\noi so that by (\ref{eqn_4_6}) $A = \vol(\Delta) I_n$ which is symmetric and positive definite and hence the set of centroids is immobilizing. We can also see this by noting that for $G$ the corresponding matrix $\Lambda$ is $\frac{1}{n}(J - I)$ where $J$ is the $(n+1) \times (n+1)$ matrix all of whose entries are 1. $J$ has one eigenvalue of $n+1$ and $n$ eigenvalues of 0, so that $\Lambda$ has one eigenvalue of 1 and $n$  eigenvalues of $-\frac{1}{n},$ whence by Lemma \ref{lemma7}, $A$ has n eigenvalues of $\vol(\Delta).$
 
\subsection{Centred contact sets} \label{subsec:centred}

\noi A set of contact points $\mathcal{P} = \{p_0, p_1, \ldots, p_n\}$ where $p_i \in F_i$ is centred at $z \in \Bbl{R}^n$ if the normal lines $l_i: x = p_i + t \mathbf{k}_i, ~ t \in \Bbl{R}$ are concurrent, meeting in $z.$ Given $\Delta$, the set $\mathcal{Z}$ of such points $z \in \Bbl{R}^n$ may be constructed as the intersection of the $n+1$ cylinders $C_i$ where $C_i = \{l_i: p_i \in F_i\}.$ This set is a polytope in $\Bbl{R}^n$ which we show in Appendix B always includes an open set of points interior to $\Delta,$ but not necessarily the whole of the interior points of $\Delta$ (for $n = 2,$ $\Delta$ a triangle, if $\Delta$ is acute angled, then the set is a hexagon completely enclosing $\Delta$, of twice the area of $\Delta;$ if $\Delta$ is obtuse angled, the set is a parallelogram only including part of $\Delta$ and which becomes unbounded as the obtuse angle approaches $\pi$).

\noi Let $z$ be a point of concurrency interior to $\Delta.$ Then there exists $t_i > 0$ so that $p_i = z + t_i \mathbf{k}_i \in F_i$ and  for the contact set $\mathcal{P} = \{p_i\}$
$$A = \sum_{i=0}^n \mathbf{k}_i z^T + \sum_{i=0}^n t_i \mathbf{k}_i \mathbf{k}_i^T =  \sum_{i=0}^n t_i \mathbf{k}_i \mathbf{k}_i^T$$
which is symmetric and positive definite, so that $\mathcal{P}$ immobilizes $\Delta.$

\medskip
\begin{rem} For $n = 2, \Delta$ a triangle, the collection of contact sets centred at $z$ whether or not $z$ lies within $\Delta$, is precisely the collection of all immobilizing contact sets \cite{CSU}. This is not true for $n \geq 3.$ From the results of the appendices, the space of centred immobilizing contact sets has dimension $n$ whereas the space of all immobilizing contact sets has dimension $\frac{1}{2}n(n+1) - 1.$ These are equal when $n = 2$ but not otherwise. \end{rem}

\subsection{Displacements} \label{subsec:displacements}
A comparatively complete description of immobilizing contact sets  of $\Delta$ can be given by examining, instead of the contact set  itself, displacements $\Delta p_i$ of the points $p_i$ within the faces $F_i$ from one immobilizing contact set to another. Leaving aside the condition of almost positive definiteness for the time being, the conditions $p_i$ lying in the hyperplane containing $F_i$ and $A$ symmetric are:
$$
\mathbf{k}_i \cdot  p_i = \mathbf{k}_i \cdot v_j, ~ j \neq i, ~ i = 0, 1, \ldots, n; ~~~ \sum_{i=0}^n \mathbf{k}_i \wedge p_i = \mathbf{0}, 
$$
\noi the latter merely stating that the antisymmetric part of $A$ vanishes. Then a displacement $\Delta p_i$ within $F_i$ maintains the symmetry of $A$ only if $\Delta p_i$ satisfy
\begin{equation}
\mathbf{k}_i \cdot \Delta p_i = 0, ~~~ i = 0, 1, \ldots, n, ~~~~~ \sum_{i=0}^n \mathbf{k}_i \wedge \Delta p_i = \mathbf{0}. \label{eqn_5_3_1} 
\end{equation}

\noi While the system (\ref{eqn_5_3_1}) is necessary to keep each  displaced point $p_i$ in the face $F_i$ and to maintain symmetry of $A,$ it is not sufficient, because it places no constraint on the size of the displacements $\Delta p_i:$ the system (\ref{eqn_5_3_1}) achieves the more modest objective of keeping each point $p_i$ in the hyperplane containing the face $F_i$ (equivalently, the displacement $\Delta p_i$ is parallel to $F_i$) together with maintaining symmetry of $A.$ 

\noi We now exhibit a basis for the space $\Omega$ of displacements $\Delta P =$ $[\Delta p_0,$ $\Delta p_1,$ $\ldots, \Delta p_n]$ satisfying the system (\ref{eqn_5_3_1}).   There is no such displacement set where all but one of the $\Delta p_i$ are zero, since, if for example, $\Delta p_0$ were the only non-zero displacement, we would have $\mathbf{k}_0 \cdot \Delta p_0 = 0$ and $\mathbf{k}_0 \wedge p_0 = 0$ which requires $\Delta p_0$ to be both perpendicular and parallel to $\mathbf{k}_0.$\\
There are, however, displacements where all but two of the $\Delta  p_i$ are zero.
 If the non-zero displacements are $\Delta p_i$ and $\Delta p_j,$ then (\ref{eqn_5_3_1}) reduces to 
 \begin{equation}
 \mathbf{k}_i \cdot \Delta p_i = 0, \mathbf{k}_j \cdot \Delta p_j = 0, \label{eqn_5_3_2} 
 \end{equation}

\begin{equation}
\mathbf{k}_i \wedge \Delta p_i + \mathbf{k}_j \wedge \Delta p_j = \mathbf{0}. \label{eqn_5_3_3} 
\end{equation}
\noi Wedging (\ref{eqn_5_3_3}) with $\mathbf{k}_i$ gives $\Delta p_j \wedge (\mathbf{k}_i \wedge \mathbf{k}_j) = \mathbf{0}$ so that $\Delta p_j$ is a linear combination of $\mathbf{k}_i$ and $\mathbf{k}_j;$ likewise $\Delta p_i.$ Then using (\ref{eqn_5_3_2}), it follows that
\begin{equation}
\Delta p_i = t \mathbf{k}_{ij} ~~\hbox{and}~~ \Delta p_j = t \mathbf{k}_{ji}, \label{eqn_5_3_4}
\end{equation}
where

\begin{equation}
\mathbf{k}_{ij} = \mathbf{k}_j - \frac{(\mathbf{k}_i \cdot \mathbf{k}_j)}{|\mathbf{k}_i|^2} \mathbf{k}_i \label{eqn_5_3_5}
\end{equation}
is the projection of $\mathbf{k}_j$  perpendicular to $\mathbf{k}_i$ and $t$ is a parameter. We define a collection of displacement sets $\Delta P_{ij}, ~0 \leq i < j \leq n$ by

\begin{equation} \Delta P_{ij} = [\Delta p_0, \Delta p_1, \ldots, \Delta p_n]; 
 ~ \Delta p_i = \mathbf{k}_{ij};~ \Delta p_j = \mathbf{k}_{ji}: ~~\\
\Delta p_m = \mathbf{0} ~\hbox{for all}~ m \neq i, j. \label{eqn_new_e1} \end{equation}

\noi We have thus shown that the displacements $\Delta P = [\Delta p_0, \Delta p_1, \ldots, \Delta p_n]$ with each $\Delta p_k$ parallel to $F_k,$ which also maintains symmetry of $A,$ and for which only $\Delta p_i$ and $\Delta p_j, i \neq j,$ are non-zero is necessarily given by $\Delta P = t \Delta P_{ij}$ for some $t$. In the displacement $\Delta P_{ij}$, the points $p_i$ and $p_j$ move in a dependent way in directions parallel to the orthogonal projection of $\mathbf{k}_j$ onto $F_i$ and the orthogonal projection of $\mathbf{k}_i$ onto $F_j$ respectively; since these directions are both linear combinations of $\mathbf{k}_i$ and $\mathbf{k}_j$, they are both orthogonal to $F_{ij}$, the $n-2$ dimensional face in which $F_i$ meets $F_j$.

\medskip
\noi In Appendix C, we show that the collection  $\mathcal{C}$ of $\Delta P_{ij}, 0 \leq i < j \leq n$ spans the space $\Omega$ of displacements which satisfy (\ref{eqn_5_3_1}), but they are not independent (see Appendix D).  

\medskip
\noi We can see the effect of the displacement set $t\Delta P_{01}$ on the immobilizing set of centroids $\mathcal{G}$. Since each $g_i$ is strictly interior to $F_i,$ for sufficiently small $|t|,$ the displaced points remain in their faces and, by construction, $A$ remains symmetric, so we examine almost positive definiteness. For the displaced set we have
$$[p_0, p_1, \ldots, p_n] = [g_0, g_1, \ldots, g_n] + t[\mathbf{k}_{01}, \mathbf{k}_{10}, 0, \ldots, 0]$$
so that 
$$A = \vol(\Delta) I_n + t B_{01} $$
where $B_{01}$ is the $n \times n$-matrix
$$B_{01} = \mathbf{k}_0\mathbf{k}^T_{01} + \mathbf{k}_1\mathbf{k}^T_{10}$$
Now take a basis of $\Bbl{R}^n$ given by $\{\mathbf{k}_0, \mathbf{k}_1, \mathbf{u}_2, \ldots, \mathbf{u}_{n-1} \}$ where each $\mathbf{u}_i$ is orthogonal to span $\{\mathbf{k}_0, \mathbf{k}_1 \}$. Then
\begin{eqnarray*}
B_{01} \mathbf{k}_0 & =  \mathbf{k}_1  |\mathbf{k}_{10}|^2 \\
B_{01} \mathbf{k}_1 & =  \mathbf{k}_0  |\mathbf{k}_{01}|^2 \\
B_{01} u_j & =  \mathbf{0}, ~j = 2, \ldots, n \\
\end{eqnarray*}
so that $B_{01}$ is similar to a matrix $(b_{ij})$ where $b_{10} = |\mathbf{k}_{10}|^2, b_{01} = |\mathbf{k}_{01}|^2$ and all other entries are zero. Thus $B_{01}$ has eigenvalues $\pm |\mathbf{k}_{01}| |\mathbf{k}_{10}|$ and all other eigenvalues are 0, so that $A$ has eigenvalues $\vol(\Delta) \pm t |\mathbf{k}_{01}| |\mathbf{k}_{10}|$ and the remaining $n-2$ eigenvalues are all $\vol(\Delta),$ hence for the displaced contact set, the smallest sum of any pair of eigenvalues is 
$$2\vol(\Delta) -  |t| |\mathbf{k}_{01}| |\mathbf{k}_{10}|$$
so that any displacement $t \Delta P_{01}, t \neq 0,$ reduces the minimal sum of a pair of eigenvalues and hence for $|t|$ sufficiently large the almost positive definite condition fails.

\medskip
\noi We extend this to an arbitrary non-zero displacement  $\Delta P \!= \! \sum_{0 \leq i < j \leq n} \! t_{ij} \Delta P_{ij}$ as follows. Following from above, we now have
$$ A = \vol(\Delta)I_n + \sum_{0 \leq i < j \leq n} t_{ij} B_{ij}$$
where
$$B_{ij} = \mathbf{k}_i \mathbf{k}^T_{ij} + \mathbf{k}_j\mathbf{k}^T_{ji} .$$
Now from before $B_{01}$ is similar to $(b_{ij})$ which is trace-free. Thus $B_{01}$ itself and similarly all $B_{ij}$ are trace-free. Hence the sum of the eigenvalues  $\lambda_1, \lambda_2, \ldots, \lambda_n$ of $A$ satisfies
$$\lambda_1 + \lambda_2 + \cdots + \lambda_n = \tr(A) = \tr(\vol(\Delta)I_n + \sum t_{ij}B_{ij}) = n \vol(\Delta). $$

\noi Now $\lambda_1 = \lambda_2 = \cdots = \lambda_n = \vol(\Delta)$ if and only if $A = \vol(\Delta)I_n$ (only if requires us to use symmetry of $A$ so that $A$ is diagonalizable). So for any non-zero perturbation $\sum t_{ij} B_{ij}$ of $A$, since the sum of eigenvalues is fixed, at least one eigenvalue increases strictly and at least one eigenvalue decreases strictly, so that the smallest sum of any pair of eigenvalues of $A$ has to decrease for any non-zero displacement $\sum t_{ij} \Delta P_{ij}$ from $\mathcal{G}.$ Because the eigenvalues of $A$ change continuously with $\Delta P,$ using the dimensionality results in the appendices, we can summarise the results as follows:

\be  \item in the $\frac{1}{2}n(n+1) - 1$ dimensional space of contact sets $[p_0, p_1, \ldots, p_n]$ for which $A$ is symmetric and where each $p_i$ lies in the $(n-1)-$dimensional hyperplane containing $F_i,$ there is an open neighbourhood $N$ of $\mathcal{G},$ the centroids, which immobilizes $\Delta;$
 \item the contact set $\mathcal{G}$ is optimal in that any displacement within the  neighbourhood $N$ causes the smallest sum of a pair of eigenvalues of $A$ to decrease from $2 \vol(\Delta),$ and if such a displacement is large enough to cause this sum to vanish, then immobilization is lost.

\ee


\noi
{\bf Acknowledgement}
The  authors wish to thank Elmer Rees for his generous advice towards this paper.

\newpage
\section{Appendix A}

 \noi  We give an example of a simplex $\Delta$ in $\Bbl{R}^4$ and points $p_0, p_1,  \ldots, p_4$ on $\Delta$ satifying the symmetry condition but failing the almost positive definiteness condition. Consider the 4-simplex having vertices $v_0 =
(-\frac{5}{12},-1,0,-3)$, $v_1 = (-\frac{83}{36}, 0,0,1)$,
$v_2 = (1,1,0,-3)$, $v_3 = (\frac{35}{18},0,-1,1)$ and $v_4 =
(\frac{35}{18},0,1,1)$. The standard outward normal vectors of the simplex
are $\mathbf{k}_0 =
(0, 34, 0, \frac{17}{2})$, $\mathbf{k}_1 =
(16, -\frac{34}{3}, 0, -\frac{119}{18})$, 
$\mathbf{k}_2 = (0, -34, 0, \frac{17}{2})$, $\mathbf{k}_3 =
(-8, \frac{17}{3}, 34, -\frac{187}{36})$ and \\
$\mathbf{k}_4 = (-8, \frac{17}{3}, -34, -\frac{187}{36})$. The points 
\begin{eqnarray}
p_0 &= \frac{3}{10}v_1 + \frac{2}{5}v_2 + \frac{3}{20}v_3 +
\frac{3}{20}v_4\\
p_1 &=  \frac{1}{10}v_0 + \frac{1}{10}v_2 + \frac{2}{5}v_3 +
\frac{2}{5}v_4\\
p_2 &=   \frac{2}{5}v_0 + \frac{2}{5}v_1 + \frac{1}{10}v_3 +
\frac{1}{10}v_4\\
p_3 &=   \frac{1}{10}v_0 + \frac{7}{10}v_1 + \frac{1}{10}v_2 +
\frac{1}{10}v_4\\
p_4 &=   \frac{1}{10}v_0 + \frac{7}{10}v_1 + \frac{1}{10}v_2 +
\frac{1}{10}v_3
\end{eqnarray} 
\noi are interior to their faces and satisfy the symmetry condition since 
$$\sum_{i=0}^{4} \mathbf{k}_ip_i^T = \left[ \begin{array}{cccc}
\frac{238}{5} & 0 & 0 & 0\\ 0 & \frac{136}{5} & 0 & 0\\ 0 & 0 &
\frac{34}{5} & 0 \\ 0 & 0 & 0 & \frac{-68}{5} \end{array} \right]. $$
\noi However, a pair of eigenvalues of this matrix has a negative sum.

\newpage

\section{Appendix B}
\noi We demonstrate that the set $\mathcal{Z}$ of points of concurrency of centred contact sets $\mathcal{P}$ defined in \S~\ref{geometric} contains an open subset of $\Bbl{R}^n$ which lies in the interior of $\Delta$.

\noi Through any $z \in \Bbl{R}^n$ pass $n+1$ lines $\ell_j, j = 0, 1,  \ldots, n,$ where $\ell_j$ has direction $\mathbf{k}_j$, so meets orthogonally $\pi_j$, the hyperplane containing $F_j.$ The line $\ell_j$ meets $\pi_j$ in a point $z_j$ satisfying

\begin{equation} z_j = z + t_j \mathbf{k}_j \label{append1}
\end{equation}

and
\begin{equation} \mathbf{k}_j \cdot z_j = \mathbf{k}_j \cdot v_i ~~\hbox{for any}~~i \neq j \label{app_a2} \end{equation} 

\noi The requirements that each $z_j$ lies in $F_j$ and the aim to show that there exists such points $z$ within $\Delta$ are best encoded by expressing the points $z_i$ and $z$ as linear combinations of the vertices $[v_0, v_1, \ldots, v_n].$ Thus we write

\begin{equation} z_j = \sum_{i=0}^n \lambda_{ij} v_i \label{app_a3} \end{equation}
where $(\lambda_{ij})$ is a stochastic matrix just as in \S~\ref{matrix} and we seek

\begin{equation} z = \sum_{i=0}^n \mu_i v_i \label{app_a4} \end{equation}

\noi where $\sum_{i=0}^n \mu_i = 1$ and $\mu_i > 0$ to ensure that $z$ lies within $\Delta.$ The combination of $v_i$ and $\mathbf{k}_j$ in the system (\ref{append1})  - (\ref{app_a4})
 suggest that it is advantageous to use the machinery developed in \S~\ref{matrix} so we again extend into $\Bbl{R}^{n+1}$ by writing $\bar{z} = (1, z)$ and $\bar{z}_i = (1, z_i).$ Then (\ref{app_a3}) and (\ref{app_a4}) give
\begin{equation}
[\bar{z}_0, \bar{z}_1, \ldots, \bar{z}_n] = V \Lambda \label{app_a5}
\end{equation}
and 
\begin{equation}
\bar{z} = V \mu,
\end{equation}
\noi where $\mu$ is the $(n+1)$ column $(\mu_0, \mu_1, \ldots, \mu_n).$ Then from the equation $K^TV = - n \vol(\Delta) I,$ we have
\begin{equation}
-n \vol(\Delta) \lambda_{ij} = \bar{\mathbf{k}}_i^T \bar{z}_j \label{app_a7}
\end{equation}
and
\begin{equation}
-n \vol(\Delta) \mu_i = \bar{\mathbf{k}}_i^T \bar{z}. \label{app_a8}
\end{equation}

\noi From \S~3 equation (\ref{eqn_2_4}), for all $i \neq j$,
\begin{equation}
\label{35A} \mathbf{k}_j \cdot v_i = - \kappa_j,
\end{equation}
so by (\ref{append1}) dotted with $\mathbf{k}_j$, (\ref{app_a2}) and (\ref{35A}) we have
$$ \mathbf{k}_j \cdot z + t_j |\mathbf{k}_j|^2 = - \kappa_j $$
\noi from which it follows
\begin{equation}
t_j = - \frac{\mathbf{k}_j \cdot z + \kappa_j}{|\mathbf{k}_j|^2} = - \frac{\bar{\mathbf{k}}_j^T \bar{z}}{|\mathbf{k}_j|^2} =  + \frac{ n \vol(\Delta) \mu_j}{|\mathbf{k}_j|^2} \label{app_a9}
\end{equation}
\noi by means of (\ref{app_a8}) and where we recall that $\bar{\mathbf{k}}_j = (\kappa_j, \mathbf{k}_j).$ Hence (\ref{append1}) becomes
\begin{equation}
z_j = z + \frac{n \vol(\Delta)}{|\mathbf{k}_j|^2} \mu_j \mathbf{k}_j, \label{app_a10}
\end{equation}
which, given $z \in \Bbl{R}^n,$ locates the  $z_j \in \pi_j. $

\medskip
\noi We want to extend equation (\ref{app_a10}) into $\Bbl{R}^{n+1}$ so that we can use (\ref{app_a7}) and (\ref{app_a8}) to find a relation between the coefficients $\mu$ and $\Lambda.$ As the zero components of both $\bar{z}_j$ and $\bar{z}$ are both unity we cannot simply replace $\mathbf{k}_j$ by $\bar{\mathbf{k}}_j$ in (\ref{app_a10}); we will need to extend $\mathbf{k}_j$ to $\Bbl{R}^{n+1}$ in such a way that the zero component vanishes. For this we use $\sum_l \mathbf{k}_l = \mathbf{0}$ and replace $\mathbf{k}_j$ by
$$\mathbf{k}_j + \rho_j \sum_l \mathbf{k}_l$$
\noi and then choose $\rho_j$ so that the 0 component in the corresponding expression
$$\bar{\mathbf{k}}_j + \rho_j \sum_l \bar{\mathbf{k}}_l$$
\noi vanishes. Thus we require
$$0 = \kappa_j + \rho_j \sum_l \kappa_l = \kappa_j + \rho_j(-n \vol(\Delta)),$$
\noi using (\ref{eqn_2_5}) and now (\ref{app_a10}) may be extended into $\Bbl{R}^{n+1}$ as 

\begin{equation}
\bar{z}_j = \ol{z} + \frac{n \vol(\Delta)}{|\mathbf{k}_j|^2} \mu_j \left( \bar{\mathbf{k}}_j + \frac{\kappa_j}{n \vol(\Delta)} \sum_l \bar{\mathbf{k}}_l \right) \label{app_a11}
\end{equation}
\noi Now multiplying by $\bar{\mathbf{k}}_i^T$ and using (\ref{app_a7}) and (\ref{app_a8}), there follows
\begin{equation}
 -n \vol(\Delta) \lambda_{ij} = - n \vol(\Delta) \mu_i + \frac{n \vol(\Delta)}{|\mathbf{k}_j|^2} \mu_j \left( \bar{\mathbf{k}}_i^T \bar{\mathbf{k}}_j + \frac{\kappa_j}{n \vol(\Delta)} \bar{\mathbf{k}}_i^T \left(\sum_l \bar{\mathbf{k}}_l\right) \right). \label{app_a12}
\end{equation}

\noi Now for the final term above  $\bar{\mathbf{k}}_i^T  \bar{\mathbf{k}}_j = \kappa_i \kappa_j +  \mathbf{k}_i \cdot \mathbf{k}_j$ while  $\bar{\mathbf{k}}_i^T \sum_l  \bar{\mathbf{k}}_l =  \bar{\mathbf{k}_i} \cdot \left( \sum_l \kappa_l, \mathbf{0} \right)$ since $\sum_l \mathbf{k}_l = \mathbf{0},$ and this equals
$\kappa_i \sum_l \kappa_l = -n \vol(\Delta) \kappa_l$ so that the final bracket in (\ref{app_a12}) just reduces to $\mathbf{k}_i \cdot \mathbf{k}_j$ and thus (\ref{app_a12}) becomes

\begin{equation}
\lambda_{ij} = \mu_i - \frac{\mu_j}{|\mathbf{k}_j|^2} \mathbf{k}_i \cdot \mathbf{k}_j. \label{app_a13}
\end{equation}
\noi We note that $\lambda_{jj} = 0$ as required and for each $j, \sum_{i=0}^n \mu_i = 1$ gives $\sum_{i=0}^n \lambda_{ij} = 1$ (where we use $\sum_i \mathbf{k}_i = \mathbf{0}).$ For a point $z$ interior to $\Delta$ to yield the corresponding point $z_i$ interior to $F_i$ we need to find $\mu_i > 0, i = 0, 1, \ldots, n$ so that $\lambda_{ij} > 0$ for all $0 \leq i \neq j \leq n.$ A particular solution to this is
\begin{equation}
\mu_i = \frac{|\mathbf{k}_i|}{\sum_{l=0}^n |\mathbf{k}_l|} \label{eqn_app14}
\end{equation}
\noi since then
$$\lambda_{ij} = \frac{|\mathbf{k}_j|\left(|\mathbf{k}_j||\mathbf{k}_i|  - \mathbf{k}_i \cdot \mathbf{k}_j\right)}{|\mathbf{k}_j|^2 \sum_{l=0}^n |\mathbf{k}_l|}$$
\noi which is positive for all $i \neq j.$ For the point $z$ corresponding to the solution (\ref{eqn_app14}), by continuity, there is a full neighbourhood $N$ of $z$ interior to $\Delta$ so that each $x \in N$ projects along $\mathbf{k}_i$ to a point $x_i$ interior to $F_i$ and thus yields an immobilizing set of contact points.

\medskip
\noi We note that for $n = 2,$ for the centroid immobilizing set where $\lambda_{ij} = \frac{1}{2}$ for $0 \leq i \neq j \leq 2,$ equation (\ref{app_a13}) has  a solution given by
$$\mu_0 = \frac{-\frac{1}{2} \mathbf{k}^2_0 \mathbf{k}_i \cdot \mathbf{k}_2}{\Delta^2}$$
and so on cyclically, where $\Delta^2$ is $|\mathbf{k}_i \times \mathbf{k}_j|^2$ for any pair $i \neq j.$ This shows that for the triangle, the centroid contact set is a centred contact set (which is just a complicated way of saying that the perpendicular bisectors of the sides of a triangle are concurrent). There is no corresponding result for $n \geq 3$ so that, in genenal, for $n \geq 3$ the centroid contact set is not a centred contact set.

\newpage
\section{Appendix C}
\noi We demonstrate that the collection $\mathcal{C}$ of displacements $\Delta P_{ij}, 0 \leq i < j \leq n$ defined by (\ref{eqn_5_3_5}) and (\ref{eqn_new_e1}) in \S~\ref{geometric} spans the space $\Omega$ of displacements
$\Delta P = [\Delta p_0, \Delta p_1, \ldots, \Delta p_n]$ with each $\Delta p_i$ parallel to $F_i$ and which also maintains symmetry of $A$; that is those $\Delta P$ which satisfy the system (\ref{eqn_5_3_1}). 

\medskip
\noi We firstly observe that for each $i = 0, 1, \ldots, n$ the $\mathbf{k}_{ij}, j \neq i$ span the space of directions parallel to $F_i$: a vector $v$ lies in this space if and only if $\mathbf{k}_i \cdot v = 0.$ Now $v = \sum_j v_{ij} \mathbf{k}_j $ since the $\mathbf{k}_j$ span $\Bbl{R}^n$ so that
$\sum_j v_{ij}(\mathbf{k}_i \cdot \mathbf{k}_j) = 0$ whence
$$v_{ii} = - \sum_{j \neq i} \frac{\mathbf{k}_i \cdot \mathbf{k}_j}{|\mathbf{k}_i|^2} v_{ij}$$ giving
$v = \sum_{j \neq i} v_{ij} \mathbf{k}_{ij}.$ For each $i$, the set $\{ \mathbf{k}_{ij}: j \neq i\}$ is dependent: from $\sum_{j=0}^n \mathbf{k}_j = \mathbf{0}$ there follows the unique dependency
\begin{equation}
\sum_{j \neq i} \mathbf{k}_{ij} = \mathbf{0}. \label{eqn_b3}
\end{equation}
\noi Thus each $\Delta P$ having each $\Delta p_i$ parallel to $F_i$ may be written (non-uniquely) as \\ $\Delta P = [\Delta p_0, \Delta p_1, \ldots, \Delta p_n]$ where
\begin{equation}
\Delta p_i = \sum_{j=0}^n c_{ij} \mathbf{k}_{ij} \label{eqn_b4}
\end{equation}
in which $c_{ii} = 0,$ and this ensures that the parallel condition $\mathbf{k}_i \cdot \Delta p_i = 0$ is satisfied.

\medskip
\noi We now exploit the lack of uniqueness in the representation (\ref{eqn_b4}). For $i = 0, 1, \ldots, n,$ in view of (\ref{eqn_b3}) we have
$$\Delta p_i = \sum_{j=0}^n c_{ij} \mathbf{k}_{ij} + \mu_i \sum_{j=0, j \neq i}^n \mathbf{k}_{ij}$$ and we choose $\mu_i$ so that
$$c_{i0} + \mu_i = c_{0i}~~\hbox{for}~ i = 0, 1, \ldots, n,$$
(so that $\mu_0 = 0$) and now define $c'_{ij}, ~ 0 \leq i, j \leq n$ by 
\begin{eqnarray}  c'_{0i} &=   c_{0i},~ i = 0, 1, \ldots, n; \\
c'_{i0} &=   c_{i0} + \mu_i,~ i = 0, 1, \ldots, n; \\
c'_{ij} &=   c_{ij} + \mu_i,~ 1 \leq i \neq j \leq n; \\
c'_{ii} &=   0,~ 0 \leq i \leq n. \\
\end{eqnarray} 
\noi Then for the given $\Delta P$ we now have
\begin{equation}
\Delta p_i = \sum_{j=0}^n c'_{ij} \mathbf{k}_{ij} \label{eqn_b5}
\end{equation}
where $c'_{ii} = 0$ and $c'_{i0} = c'_{0i}$ for $i = 0, 1, \ldots, n.$ We now study the effect of the symmetry condition 
\begin{equation}
\sum_{i=0}^n \mathbf{k}_i \wedge \Delta p_i = 0 \label{eqn_b6}
\end{equation}
on the representation (\ref{eqn_b5}) and observe  from (\ref{eqn_5_3_5}) that
\begin{equation}
\mathbf{k}_i \wedge \mathbf{k}_{ij} = \mathbf{k}_i \wedge \mathbf{k}_j. \label{eqn_b7}
\end{equation}
Hence substituting (\ref{eqn_b5}) into (\ref{eqn_b6}) and using (\ref{eqn_b7}) there follows

\begin{eqnarray}  0 &=  \sum_{i=0}^n \mathbf{k}_i \wedge \sum_{j=0}^n c'_{ij} \mathbf{k}_{ij} = \sum_{i, j = 0}^n c'_{ij} \mathbf{k}_i \wedge \mathbf{k}_j \\
&=   \sum_{j=0}^n c'_{0j} \mathbf{k}_{0} \wedge \mathbf{k}_j + \sum_{i = 0}^n c'_{i0} \mathbf{k}_i \wedge \mathbf{k}_0  + \sum_{i, j = 1}^n c'_{ij} \mathbf{k}_i \wedge \mathbf{k}_j \\
&=  \sum_{i, j = 1}^n c'_{ij} \mathbf{k}_i \wedge \mathbf{k}_j, 
\end{eqnarray} 
the terms with $i$ or $j$ being zero cancelling because of the skew symmetry of $\mathbf{k}_0 \wedge \mathbf{k}_j$ and the normalisation $c'_{0j} = c'_{j0}.$ But since the wedge products 
$\mathbf{k}_i \wedge \mathbf{k}_j$  are independent for $1 \leq i < j \leq n$ it follows that $c'_{ij} = c'_{ji}$ for all $1 \leq i < j \leq n,$ whence $c'_{ij} = c'_{ji}$ for $0 \leq i < j \leq n$ together with $c'_{ii} = 0, ~0 \leq i \leq n.$ Thus

\begin{eqnarray}  \Delta P = [\Delta p_0, \Delta p_1, \ldots, \Delta p_n] &=  \left[\sum c'_{0j} \mathbf{k}_{0j}, \sum c'_{1j} \mathbf{k}_{1j}, \ldots, \sum c'_{nj} \mathbf{k}_{nj} \right]\\
&=  \sum_{0 \leq i < j \leq n} c'_{ij} \Delta P_{ij} 
\end{eqnarray} 
\noi showing that the $\Delta P_{ij},~ 0 \leq i < j \leq n$ span $\Omega.$

\newpage
\section{Appendix D}
We show that the displacements $\Delta P_{ij}, 0 \leq i < j \leq n$ which satisfy equation (\ref{eqn_5_3_1}) are  not independent. If $\sum_{0 \leq i < j \leq n} c_{ij} \Delta P_{ij}$ vanishes, then each $\Delta p_m = \mathbf{0}.$ Noting that $\Delta P_{ij}$ only has non-zero components in the $i$th and $j$th locations, for the linear combination $\sum_{0 \leq i < j \leq n} c_{ij} \Delta P_{ij},$ we have for $0 \leq m \leq n$
$$\mathbf{0} = \Delta p_m = \sum_{i=0}^{m-1} c_{im} \mathbf{k}_{mi} + \sum_{j=m+1}^n c_{mj} \mathbf{k}_{mj},$$
where  for $m = 0,$ the first sum in the right hand side vanishes and for $m = n,$ the second sum vanishes. Hence  by (\ref{eqn_5_3_5}) for $0 \leq m \leq n,$

$$\sum_{i = 0}^{m-1} c_{im} \mathbf{k}_i + \sum_{j=m+1}^n c_{mj} \mathbf{k}_j - \left( \sum_{i=0}^{m-1} c_{im} \frac{\mathbf{k}_i \cdot \mathbf{k}_m}{|\mathbf{k}_m|^2} +  \sum_{j=m+1}^{n} c_{mj} \frac{\mathbf{k}_j \cdot \mathbf{k}_m}{|\mathbf{k}_m|^2} \right) \mathbf{k}_m = \mathbf{0}.$$

\noi But the only dependency between $\mathbf{k}_0, \mathbf{k}_1, \ldots, \mathbf{k}_n$ is $ \sum_{j=0}^n \mathbf{k}_j = \mathbf{0}$ from which it follows that 
$$c_{0m} = c_{1m} = \cdots = c_{(m-1)m} = c_{m(m+1)} = \cdots = c_{mn} .$$
\noi Applying this for each $m = 0, 1, \ldots, n$ gives that all the $c_{ij}, 0 \leq i < j \leq n$ must be equal. Thus there is precisely one dependency between the $\Delta P_{ij}$ and the space of displacement sets maintaining symmetry of $A$ is $\frac{1}{2} n (n+1)-1$ dimensional.

\end{document}